\documentclass[12pt]{article}

\title{Fires on trees} 
\author{Jean Bertoin\thanks{Laboratoire de Probabilit\'es et Mod\`eles Al\'eatoires, 
UPMC, 4 Place Jussieu, 75252 Paris Cedex 05; France
Email: jean.bertoin@upmc.fr} 
 }

\date{}
\usepackage{graphicx}
\usepackage{amsfonts,amsmath,amssymb,bbm}
\usepackage{setspace} 
\def\proof{\noindent{\bf Proof:}\hskip10pt}        
\def\QED{\hfill $\Box$}

\font\tenmath=msbm10 scaled 1200
\font\sevenmath=msbm7 scaled 1200
\font\Fivemath=msbm5 scaled 1200
\newfam\mathfam \textfont\mathfam=\tenmath
\scriptfont\mathfam=\sevenmath \scriptscriptfont\mathfam=\Fivemath

\def \\ { \cr }

\def \1{1 \mkern -6mu 1} 
\def\N{\mathbb{N}}
\def\E{\mathbb{E}}
\def\P{\mathbb{P}}
 
\def\f{\mathcal{F}_n}

\def \e{{\rm e}}
\def \d{{\rm d}}

\def \t{{\bf t}_n}

\newtheorem{theorem}{Theorem}

\newtheorem{proposition}{Proposition}
\newtheorem{lemma}{Lemma}
\newtheorem{corollary}{Corollary}
\begin{document}

\maketitle

\begin{abstract} 
We consider random dynamics on the edges of  a uniform Cayley tree with $n$ vertices, in which edges are either inflammable, fireproof, or burnt.
Every inflammable edge is replaced by a fireproof edge at unit rate, while
fires start at smaller rate $n^{-\alpha}$ on each inflammable edge, then propagate through the neighboring inflammable edges and are only stopped at fireproof edges. A vertex is called fireproof when all its adjacent edges are fireproof. We show that as $n\to \infty$, the density of fireproof vertices converges to $1$ when $\alpha>1/2$, to $0$ when $\alpha<1/2$, and to some non-degenerate random variable when $\alpha=1/2$. We further study the connectivity of the fireproof forest, in particular the existence of a giant component. 
\end{abstract}


\begin{section}{Introduction }

Since their introduction  by Drossel and Schwabl \cite{DS}, forest fire models have generated a lot of attention in the literature in statistical physics. In this work, we will consider the following variation. Imagine that the edges of some finite graph can be in either of three states: inflammable, fireproof, or burnt. We suppose that initially every edge is inflammable,
and that the system evolves randomly as follows when time passes. The only transitions are from
inflammable to fireproof or burnt.  An inflammable edge becomes fireproof  at unit rate and independently of the other edges.  On the other hand, an inflammable edge can also be set on fire at rate $r>0$, called the firing rate. In that case the fire propagates to the neighboring inflammable edges and the propagation only ceases at fireproof edges. After some time, all edges are either burnt or fireproof.

Our interest in the model above has been triggered by a recent work of R\'ath \cite{Rath} (see also \cite{MR}) on so-called mean field frozen percolation. More precisely, R\'ath considers a variant of the Erd\H os-R\'enyi random graph on a large set of vertices in which lightings strike vertices at some small rate.
A lighting burns the entire connected component of the vertex that is hit, and that burnt component is then removed from the graph.  Therefore  the present model
may be viewed as  dual to mean field frozen percolation, without creation of edges but with the introduction of barriers that stop the propagation of fires.

In this work, we assume that the graph on which the fire dynamics occur is a uniform Cayley tree of size $n$, $\t$, where $n$ is some large integer.
This means that $\t$ is picked uniformly at  random amongst the $n^{n-2}$ different trees on a set of $n$ labeled vertices, say $[n]=\{1,\ldots, n\}$.
This simple choice is motivated by the fact that clusters in the random graph model  before the phase transition resemble uniform random Cayley trees.
We also suppose that the firing rate is $r=n^{-\alpha}$ where $\alpha >0$ is some parameter.

A vertex is called fireproof if all its adjacent edges are fireproof in the final state of the system, and burnt otherwise. In other words, a vertex is burnt  as soon as one of its adjacent edge has burned. We are interested in the subgraph formed by fireproof vertices, and fireproof edges between fireproof vertices. 
The fireproof edges having one or two burnt extremities and burnt edges are discarded. 
Our first result shows that the system exhibits a phase transition at $\alpha=1/2$. As $n\to\infty$, the density of fireproof vertices converges in distribution to a non-degenerate random variable when $\alpha=1/2$,  to $1$ for $\alpha >1/2$ and to $0$ for $\alpha<1/2$.
We further study the connectivity in the fireproof forest. We shall prove that
in the super critical case $\alpha >1/2$, there exists a giant tree-component of size $\sim n$ with high probability when $n\to\infty$. In the critical case $\alpha=1/2$, for every $\varepsilon >0$, with high probability there is no giant tree-component  of size at least $\varepsilon n$, although one can find
tree-components of size greater than $n^{1-\varepsilon}$. 

The rest of this paper is organized as follows. We shall start with some preliminaries on a limit theorem in distribution for  the number of random cuts which are needed to isolate a fixed number of vertices in a large uniform Cayley tree, relying  crucially on a recent work by Haas and Miermont \cite{HM2}. 
This extends earlier results of Janson \cite{Janson} and Panholzer \cite{Panholzer} and may be of independent interest.
We shall also recall well-known properties of the spinal decomposition of uniform Cayley trees, which will play an important role in the study.
Our main results about the asymptotic behavior of the density of fireproof vertices will be stated and proved in Section 3. Finally Section 4 is devoted to connectivity properties of the fireproof forest, in particular the existence of giant components, in the critical and supercritical cases.

\end{section}

\begin{section}{Preliminaries on uniform random trees}

\subsection{Random cuts and isolation of vertices}

It is sometime convenient to imagine that fireproof edges correspond to cuts on the graph which stop the propagation of fires,
and vertices are then fireproof when they have been isolated by cuts on the tree before a fire has ever reached them.
Roughly speaking, our first purpose is to estimate for every fixed integer $k\geq 1$ the number $X(n,k)$ of random cuts on a uniform Cayley tree ${\bf t}_n$ which are required to isolate 
$k$ typical vertices. 
More precisely, we  sample uniformly at random and with replacement a sequence $U_1,  \ldots, U_k$ of $k$ vertices in $[n]$, independently of ${\bf t}_n$, and consider the following stochastic algorithm.

  Imagine that we remove an edge chosen  uniformly at random given the preceding variables. This disconnects ${\bf t}_n$
 into two subtrees. If one of these two subtrees contains none of the vertices $U_1, \ldots, U_k$, then we discard it; else we keep the two subtrees.
We iterate the procedure  by removing at each step an edge uniformly at random given the current forest consisting of the subtrees that contain
at least one of the $k$ sampled vertices, and then eventually discarding the new subtree containing none of those vertices.
The algorithm terminates after $X(n,k)$ steps when all the vertices $U_1, \ldots, U_k$ have been isolated, that is when the forest reduces to the singletons
$\{U_i\}$ for $i=1, \ldots, k$.

The study of  isolation of a single point  by random cuts has been initiated by Meir and Moon \cite{MM} who have obtained the first results about the asymptotic behavior in mean of $X(n,1)$ when $n\to \infty$.
The following limit in distribution has been established in a more general setting by Janson \cite{Janson} and Panholzer \cite{Panholzer}:
$$\lim_{n\to \infty} \frac{1}{\sqrt n}X(n,1)= R\qquad \hbox{in law}$$
where $R$ denotes a Rayleigh variable, i.e.
$$\P(R\in \d s)= s \exp(-s^2/2) \d s\,, \qquad s\geq 0\,.$$
Relying on a recent paper by Haas and Miermont \cite{HM2} on scaling limits of Markov branching trees (see also Haas {\it et al.} \cite{HMPW} for a closely related work), we obtain the following extension.

\begin{lemma}\label{L1} For every integer $k\geq 1$, we have
$$\lim_{n\to \infty} \frac{1}{\sqrt n}X(n,k)= \chi(2k)\qquad \hbox{in distribution}\,,$$
where $\chi(2k)$ is a chi-variable with $2k$ degrees of freedom, viz.
$$\P(\chi(2k)\in \d x)= \frac{2^{1-k}}{(k-1)!}x^{2k-1} \exp(-x^2/2) \d x\,, \qquad x\geq 0\,.$$
\end{lemma}

\noindent {\bf Remark:} The same asymptotic result would also hold if instead, the $k$ vertices were picked  without replacement, since  sampling uniformly at random a fixed number of  vertices with replacement rather than without replacement makes no difference when the total number of vertices tends to infinity.

\proof  The proof relies crucially on the article by Haas and Miermont \cite{HM2} to whom we refer for background on some terminology and technical details which are not provided here.

We consider the logging process of some tree $\tau$ on $[n]$ by removing successively its edges in some given order. 
We call  a connected component $B\subseteq [n]$ that arises at some stage of this process a block, and view the set of blocks as another set of vertices.
Following Haas and Miermont \cite{HM1},  we represent the logging of $\tau$ as a rooted binary tree $T$ with $n$ leaves on this new set of vertices, where the root of $T$ is given by 
$[n]$ and the leaves by the singletons $\{i\}$ for $i\in[n]$. We draw an edge between a parent block $B$ and two children blocks
$B'$ and $B''$ whenever the edge-removal for the subtree on $B$ produces the two subtrees on $B'$ and $B''$, respectively. 
See Figure 1 below for an illustration.

\begin{picture}(400,280)(-10,-30)

\put(80,10){\circle{20}}
\put(40,10){\circle{20}}
\put(40,50){\circle{20}}
\put(120,50){\circle{20}}
\put(80,90){\circle{20}}
\put(160,90){\circle{20}}
\put(40,130){\circle{20}}
\put(120,130){\circle{20}}
\put(80,170){\circle{20}}

\put (40,20){\line(0,1){20}}
\put (50,10){\line(1,0){20}}
 \put (85,18){\line(1,1){26}}
 
\put (115,58){\line(-1,1){26}}
\put (75,98){\line(-1,1){26}}
 \put (85,98){\line(1,1){26}}
  \put (125,58){\line(1,1){26}}
\put (80,100){\line(0,1){60}}

\put(80,10)
{\makebox(0,0){a}}
\put(40,10)
{\makebox(0,0){h}}
\put(120,50)
{\makebox(0,0){g}}
\put(40,50)
{\makebox(0,0){f}}
\put(80,90)
{\makebox(0,0){i}}
\put(160,90)
{\makebox(0,0){d}}
\put(40,130)
{\makebox(0,0){e}}
\put(120,130)
{\makebox(0,0){c}}
\put(80,170)
{\makebox(0,0){b}}

\put(103,29){\makebox(0,0){1}}
\put(60,4){\makebox(0,0){6}}
\put(34,30){\makebox(0,0){8}}
\put(100,63){\makebox(0,0){2}}
\put(102,105){\makebox(0,0){5}}
\put(137,65){\makebox(0,0){4}}
\put(55,110){\makebox(0,0){3}}
\put(75,130){\makebox(0,0){7}}


\put(290,00){\framebox(60,15){abcdefghi} }

\put(250,50){\framebox(20,15){afh} }
\put(350,50){\framebox(40,15){bcdegi} }

\put(230,90){\framebox(10,15){a} }
\put(260,90){\framebox(15,15){fh} }  
\put(230,130){\framebox(10,15){f} }
  \put(260,130){\framebox(10,15){h} }

\put(330,90){\framebox(30,15){bcei} }
\put(400,90){\framebox(15,15){dg} }

\put(300,130){\framebox(20,15){bci} }
\put(340,130){\framebox(10,15){e} }
\put(370,130){\framebox(10,15){g} }
\put(400,130){\framebox(10,15){d} }

\put(300,170){\framebox(15,15){bi} }
\put(340,170){\framebox(10,15){c} }

\put(270,210){\framebox(10,15){b} }
\put(310,210){\framebox(10,15){i} }

\put (300,16){\line(-1,1){34}}
\put (330,16){\line(1,1){34}}

\put (260,66){\line(-1,1){24}}
\put (265,66){\line(0,1){24}}

\put (260,106){\line(-1,1){24}}
\put (265,106){\line(0,1){24}}

\put (365,66){\line(-1,1){24}}
\put (383,66){\line(1,1){24}}

\put (335,106){\line(-1,1){24}}
\put (345,106){\line(0,1){24}}

\put (400,106){\line(-1,1){24}}
\put (405,106){\line(0,1){24}}

\put (305,146){\line(0,1){24}}
\put (315,146){\line(1,1){24}}

\put (300,186){\line(-1,1){24}}
\put (313,186){\line(0,1){24}}

\end{picture}

\centerline{\bf Figure 1}
\centerline{\sl  Left: Tree $\tau$ on $9$ vertices labelled a,...,i; edges are enumerated in order of removal.}
\centerline{\sl  Right: Binary tree $T$ on the set of blocks
describing the logging of  $\tau$.}

Then select $\ell$ distinct vertices in $[n]$, say $i_1, \ldots, i_\ell$, and denote by ${\mathcal R}(T; i_1, \ldots, i_\ell)$ the smallest connected subset
of $T$ that contains the root $[n]$ and the leaves $\{i_j\}$ for $j=1, \ldots, \ell$. We call ${\mathcal R}(T; i_1, \ldots, i_\ell)$
the tree  $T$ reduced to those leaves. 
Observe that  the number of cuts in the tree $\tau$ which are needed to
isolate $i_1, \ldots, i_\ell$, in the sense of the algorithm described at the beginning of this section, 
coincides with the number of internal nodes (i.e. vertices which are not leaves) of the reduced tree ${\mathcal R}(T; i_1, \ldots, i_\ell)$. Because $T$ is binary, 
this quantity is also given by  the length (i.e. the number of edges) of the reduced tree ${\mathcal R}(T; i_1, \ldots, i_\ell)$ minus $(\ell-1)$.

We now suppose that $\tau={\bf t}_n$ is a uniform random Cayley tree on $[n]$ whose edges have been enumerated uniformly at random.
We denote by ${\bf T}_n$ the binary tree with $n$ leaves that results  as above from the logging of ${\bf t}_n$. We  view ${\bf T}_n$ endowed with the graph distance as a finite rooted metric space. Essentially, this means that we only retain the combinatorial structure of the tree ${\bf T}_n$,  forgetting  the details of vertices except for the root which is distinguishable.
We also sample $k$ points in $[n]$ uniformly at random and with replacement (so the number of different points is $\ell\leq k$), independently of the preceding, and view these points as leaves of ${\bf T}_n$.
So we have to determine the asymptotic behavior of the length of the tree ${\bf T}_n$ reduced to its branches from the root to these leaves, say ${\mathcal R}({\bf T}_n,k)$. 

The answer to this question follows from general results of Haas and Miermont \cite{HM2}. Indeed, we shall see that 
their results imply that the reduced tree ${\mathcal R}({\bf T}_n,k)$, viewed as a so-called compact rooted real tree and rescaled by a factor 
$1/\sqrt n$, converges in law in the sense of Gromov-Hausdorff  to ${\mathcal R}(k)$, the Brownian Continuous Random Tree (CRT) reduced to $k$ i.i.d. random leaves chosen 
according to its mass-measure. More generally, it would follow from Theorem 1 in \cite{HM2} that the rescaled real tree $n^{-1/2} {\bf T}_n$ converges weakly towards the Brownian CRT, but we will not need this stronger result here. 
This implies in particular that the total length of ${\mathcal R}({\bf T}_n,k)$ rescaled by $1/\sqrt n$ converges in distribution to the total length  of ${\mathcal R}(k)$.
As the latter has the chi-distribution with $k+1$ degrees of freedom (this can be seen for instance from  Lemma 21 in \cite{ACRTIII}), this yields our statement.

Hence all that we need is to check that our setting fulfills the framework of Haas and Miermont \cite{HM2}. The first point is that ${\bf T}_n$ can be viewed as a Markov branching tree with $n$ leaves, which should be plain from the work of Pitman on fragmentation of random forests; see Theorem 5 in \cite{Pi}. Next, for every positive integers  $n$ and $m$ with $n/2< m < n$, denote by $q_n(m)$ the 
probability that the first cut on ${\bf t}_n$ produces two subtrees with
sizes $m$ and $n-m$, 
respectively\footnote{The case when $n$ is even an $m=n/2$ is somewhat special due to symmetry;  it can be neglected as its probability is small when $n$ is large.}.
It is well-known and due to Pavlov \cite{Pavlov} (see also Corollary 5.8  in \cite{RFC} and Lemma \ref{L2} below) that this probability can be computed using the Borel distribution;
specifically we have
$$q_n(m)= \frac{ m^{m-1} (n-m)^{n-m-1}  (n-2)!}
{m! (n-m)! n^{n-3}}
\,,\qquad n/2< m < n\,.$$
An application of Stirling formula entails that as $n\to\infty$,
$$ n^{3/2}(1-m/n) q_n(m)\sim  \frac{1}{\sqrt{ 2\pi (m/n)^3 (1-m/n)}}\,,\qquad \hbox{ uniformly for $n/2<m<n$}\,,$$
which is the local form of Equation (3) in \cite{HM2}. It is immediate to deduce that  the basic assumption {\bf (H)} of Haas and Miermont holds
with $\gamma=1/2$, $\ell(n)\equiv 1$ and $\nu$ the measure on $\{{\bf s}=(s,1-s),1/2< s < 1\}$ given by
$$\nu(\d{\bf s})=\frac{1}{\sqrt{ 2\pi s^3(1-s)^3}}\, \d s\,.$$
Their assumption {\bf (H')} is then plain, and according to Proposition 12 in \cite{HM2}, the rescaled reduced tree $n^{-1/2}{\mathcal R}({\bf T}_n,k)$,
 viewed as a compact rooted real tree, 
converges in distribution in the sense of Gromov-Hausdorff towards the CRT denoted by ${\mathcal T}_{1/2,\nu}$ and reduced to $k$ i.i.d leaves picked according to the mass-measure on ${\mathcal T}_{1/2,\nu}$.
This completes the proof since ${\mathcal T}_{1/2,\nu}$ is the Brownian CRT, as it can be 
seen  e.g. from  pages 339-340 in \cite{Be}. \footnote{Recall that the contour process of the Brownian CRT is twice the standard Brownian excursion; which explains the extra factor $2$ in \cite{Be}.}
\QED

\subsection{Spinal decomposition} 
Our next purpose is to recall a useful decomposition of the uniform Cayley tree ${\bf t}_n$. Sample a pair of vertices $U$ and $U'$ 
uniformly at random with replacement in $[n]$ and independently of ${\bf t}_n$; $U$ should be thought of as the root. 
The oriented branch from $U$ to $U'$ is called the spine and denoted by ${\bf s}$,
its length, i.e.  the distance between $U$ and $U'$ in ${\bf t}_n$, is denoted by $\lambda_n$.
We enumerate the vertices of the spine as $V_0=U, \ldots, V_{\lambda_n}=U'$.
Removing the edges of the spine disconnects ${\bf t}_n$ into a sequence of $\lambda_n +1$ subtrees, say 
${\bf b}_0, \ldots, {\bf b}_{\lambda_n}$ which can be naturally rooted at $V_0, \ldots, V_{\lambda_n}$.
We refer  to 
the sequence $({\bf b}_0, \ldots, {\bf b}_{\lambda_n})$ as the spinal decomposition of ${\bf t}_n$
and to each ${\bf b}_i$ as a bush.

The following description of the spinal decomposition belongs to the folklore of random trees. It follows for instance from  Corollary 30  of Aldous and Pitman \cite{AP} and the observation due to Aldous \cite{ACRTIII} that  a uniform rooted Cayley tree with $n$ vertices can be viewed as a Galton-Watson tree with Poissonian offspring distribution
and conditioned to have size $n$, where  the labels are assigned to the vertices uniformly at random (and the parameter of the Poisson law irrelevant). See also Formula (104) in \cite{AP} for the law of the length of the spine. Recall that an integer valued random variable $\beta$ has the Borel$(1)$ distribution if
$$\P(\beta=k)= \frac{\e^{-k} k^{k-1}}{k!}\,,\qquad k\geq 1\,.$$

\begin{lemma} \label{L2} For every integer $n\geq 1$, we have:

\noindent {\rm (i)} The distribution of the length of the spine is given by
$$\P(\lambda_n=k)= \frac{(k+1)(n-1)!}{n^{k+1} (n-k-1)!}\,,\qquad k=0,\ldots , n-1\,.$$

\noindent {\rm (ii)} Conditionally on $\lambda_n=k$, the sizes $\beta_0, \ldots, \beta_k$ of the bushes ${\bf b}_0, \ldots, {\bf b}_{k}$ of the spinal decomposition are distributed as the sequence of $k+1$ independent Borel$(1)$ variables conditioned to have sum equal to $n$.

\noindent {\rm (iii)} Conditionally on $\lambda_n=k$ and the set of vertices ${\mathcal V}_0, \ldots , {\mathcal V}_k$ of the bushes,
 ${\bf b}_0, \ldots, {\bf b}_{k}$ are independent and each ${\bf b}_i$ has the uniform distribution on the space of rooted Cayley trees on the set of vertices
 ${\mathcal V}_i$.

\end{lemma}

\end{section}

\begin{section}{Density of fireproof vertices}
Recall the dynamics on the set of edges of ${\bf t}_n$  which has been described in the Introduction, and consider its terminal state.
For every $n\geq 1$, we write $D_n$ for the density of fireproof vertices, i.e.
$$D_n=n^{-1}{\rm Card}\{i\in[n]: \hbox{ all the edges adjacent to $i$ are fireproof}\}\,.$$

\begin{theorem} \label{T1} We have

\noindent {\rm (i)} For $\alpha >1/2$, $\lim_{n\to\infty}D_n=1$ in probability.

\noindent {\rm (ii)} For $\alpha =1/2$, 
$$\lim_{n\to\infty}D_n= D_{\infty}\qquad \hbox{in distribution}$$
where
$$\P(D_{\infty}\in \d x)=\frac{1}{\sqrt{2\pi x(1-x)^3}}\exp\left(-\frac{x}{2(1-x)}\right) \d x\,,\qquad x\in(0,1)\,.$$

\noindent {\rm (iii)} For $\alpha <1/2$, $\lim_{n\to\infty}D_n=0$ in probability.

\end{theorem}

\proof The guiding line consists in reducing the proof to Lemma \ref{L1} by viewing fireproof edges as cuts on ${\bf t}_n$. In this direction, it is 
convenient to rephrase slightly the dynamics described in the Introduction, using well-known properties of independent exponential variables.

We attach to each edge an exponential random variable with parameter $(1+n^{-\alpha})$ which we view as an alarm clock, independently of the other edges.
When the first alarm clock rings, say it is attached to the edge $e$, then we toss a coin with
$$\P({\rm Head})=1/(1+n^{-\alpha})\ \hbox{and}\  \P({\rm Tail})=n^{-\alpha}/(1+n^{-\alpha})\,.$$
If the outcome is Head, then we remove the edge $e$, thus disconnecting ${\bf t}_n$ into a pair of subtrees. If the outcome is Tail, then a fire starts at $e$
and all the edges burn instantaneously. We stress that burned edges will be kept forever. 

We iterate in an obvious way. Specifically, assume that the outcome of the first coin-tossing is Head, as otherwise we have reached the terminal state. Wait until the second alarm clock rings; suppose it is attached to the edge $e'$ which belongs to one of the two subtrees, say ${\bf t}'$. We then toss again the same coin, removing the edge $e'$ if the outcome is Head, and burning all the edges of ${\bf t}'$ if the outcome is Tail. We repeat the procedure independently each time a new alarm clock rings, agreeing that when this alarm clock is attached to an edge that has already burned, we do nothing.

A vertex $i\in[n]$ is fireproof if and only if all the edges adjacent to $i$ have been removed in the procedure above. Observe that the order on the edges that is  induced by the exponential clocks corresponds to a uniform random enumeration. Further, each time an alarm clock rings on an edge which has not yet burned, there is the same conditional probability $1/(1+n^{-\alpha})$ to remove that edge given the state of the system. Thus, if we consider arbitrary vertices 
$i_1, \ldots, i_{\ell}$, the probability that all these vertices are fireproof is given by
$$\E\left( (1+n^{-\alpha})^{-X(n; i_1, \ldots, i_{\ell})}\right)\,,$$
where $X(n; i_1, \ldots, i_{\ell})$ is the number of cuts needed to isolated each and every of the $i_j$'s when removing successively the edges of ${\bf t}_n$ in a uniform  random order and discarding the subtrees which contain none of the vertices $i_1, \ldots, i_{\ell}$.

Now pick $k$ vertices $U_1, \ldots, U_k$ uniformly at random with replacement in $[n]$ and independently of the preceding. On the one hand we have 
$$\P(U_1, \ldots , U_k \hbox{ are fireproof})=\E(D_n^k)\,,$$
where $D_n$ is the density of fireproof vertices.
On the other hand, the discussion above shows that, in the notation of Section 2.1,
$$\P(U_1, \ldots , U_k \hbox{ are fireproof})= \E\left( (1+n^{-\alpha})^{-X(n,k)}\right)\,.$$
We deduce from Lemma \ref{L1} that
$$\lim_{n\to \infty} \E(D_n^k)= 
\left\{\begin{matrix}
 1\  &\hbox{if }& \alpha >1/2\,,\\
 \E(\exp(-\chi(2k)))Ê&\hbox{if }& \alpha =1/2\,, \\
 0\  &\hbox{if }& \alpha< 1/2\,,
 \end{matrix} \right. $$
where $\chi(2k)$ is a chi-variable with $2k$ degrees of freedom. 

In particular, this proves (i) and (iii). In the critical case $\alpha = 1/2$, this also establishes   the weak convergence of $D_n$ towards some random variable $D_{\infty}$ in $[0,1]$
whose law is characterized by its entire moments, specifically
$$\E(D_{\infty}^k)=  \E(\exp(-\chi(2k)))\,,\qquad k\geq 1\,.$$
Lemma \ref{L3} below claims that the law of $D_{\infty}$ is then given as in Theorem \ref{T1}, thus completing its proof. \QED

We still have to establish the following result.

\begin{lemma} \label{L3}
Let $D_{\infty}$ be a random variable with entire moments given by 
$$\E(D_{\infty}^k)=  \E(\exp(-\chi(2k)))\,,\qquad k\geq 1\,,$$
where $\chi(2k)$ is a chi-variable with $2k$ degrees of freedom. Then
$$\P(D_{\infty}\in \d x)=\frac{1}{\sqrt{2\pi x(1-x)^3}}\exp\left(-\frac{x}{2(1-x)}\right) \d x\,,\qquad x\in(0,1)\,.$$
\end{lemma}

\proof The Laplace transform of chi-variables can be expressed in terms of Kummer's confluent hypergeometric functions; however this is not really useful to determine explicitly the law of $D_{\infty}$. We rather use the fact  that $\chi(2k)$ arises as the total length of the Brownian CRT reduced to its root and $k$ i.i.d. leaves which are picked according to the mass-measure of the CRT, as it has been pointed out in the proof of Lemma \ref{L1}. 

We now recall the fragmentation of the Brownian CRT introduced by Aldous and Pitman \cite{AP2}. Consider a Poisson point process on the skeleton of the Brownian CRT with intensity given by its length measure. We should think of the atoms of this measure as cuts on the skeleton.
Let $L_0$ be the root of the CRT and  $L_1, \ldots, L_k$ be a sequence of $k$ i.i.d leaves sampled according to the mass-measure. 
Introduce the event $\Lambda_k$ that there is no atom of the Poisson point measure on the branches connecting the root $L_0$ to the $L_j$'s. In particular, viewing $\chi(2k)$ as the length of the Brownian CRT 
reduced to $L_0, L_1, \ldots, L_k$, we have now 
$$\E(D_{\infty}^k)=  \E(\exp(-\chi(2k)))=\P(\Lambda_k).$$

The right hand side can also be expressed as $\E(Y_*^k)$, where $Y_*$ denotes the mass of the connected component of the root $L_0$ after logging the CRT at the points of the Poisson random measure. As a consequence, $D_{\infty}$ and $Y_*$ have the same law.
This establishes our claim since, according to Corollary 5 of Aldous and Pitman \cite{AP2}, the distribution of $Y_*$ is given by
 $$\P(Y_*\in \d x)=\frac{1}{\sqrt{2\pi x(1-x)^3}}\exp\left(-\frac{x}{2(1-x)}\right) \d x\,,\qquad x\in(0,1)\,.$$
\QED

\noindent{\bf Remark.} Slightly more generally, we may consider a critical fire process on a uniform Cayley tree with size $n$ with firing rate $a/\sqrt n$ for some positive constant $a$. The same arguments as in the present section show that the limiting density of fireproof vertices is given by the size of the connected component of the root of the Brownian CRT  in the fragmentation of Aldous and Pitman, where the rate of cuts on the skeleton is now $a$ times the length measure. The distribution of this limiting density
is then
$$\frac{a}{\sqrt{2\pi x(1-x)^3}}\exp\left(-\frac{a^2x}{2(1-x)}\right) \d x\,,\qquad x\in(0,1)\,.$$

\end{section}

\begin{section}{Giant components in the fireproof forest}

For every integer $n\geq 1$, we denote by $\f$ the fireproof forest which results from the dynamics on edges of the uniform Cayley tree $\t$
described in the Introduction. More precisely, the vertices of $\f$ consists of the subset of fireproof vertices, and the edges of $\f$ are the fireproof edges with fireproof end-points (recall that edges of $\t$ with at least one burnt extremity are discarded).

Our goal in this Section is to investigate some connectivity properties of $\f$ when $n$ is large. 
We are especially interested in the existence of giant components in $\f$, i.e. trees with size of order $n$. As we know from Theorem \ref{T1}
that in the subcritical case $\alpha<1/2$, the total size of $\f$ is $o(n)$ in probability when $n\to\infty$, we focus on the case $\alpha\geq 1/2$.
Observe also that for $\alpha >1$, the probability that all the $n-1$ edges of $\t$ become fireproof is high when $n$ is large, so the range of interest for $\alpha$ is $ [1/2,1]$.

In order to investigate this problem, 
we introduce  $U$ and $U'$, two independent uniform random vertices in $[n]$, independently of $\f$. 
Our main result provides an asymptotic estimate for the probability that $U$ and $U'$  
 are two connected vertices in $\f$ as $n\to\infty$.
\begin{theorem}\label{T2}
The probability that $U$ and $U'$ belong to the same tree-component
of the fireproof forest $\f$ converges when $n\to\infty$ to $1$ in the supercritical case $\alpha>1/2$,
and to $0$ in the critical case $\alpha=1/2$.
\end{theorem}

Before tackling the proof of Theorem \ref{T2}, we present an immediate consequence in terms of the existence of giant tree-components. 

\begin{corollary}\label{C1} 
Fix $\varepsilon\in(0,1)$.

\noindent {\rm (i)} In the supercritical case $\alpha>1/2$, the probability that there exists a tree-component
with size greater than $(1-\varepsilon)n$  in the fireproof forest $\f$ converges to $1$ as $n\to\infty$. 

\noindent {\rm (ii)} In the critical case $\alpha=1/2$, the probability that there exists a tree-component
with size greater than $\varepsilon n$  in the fireproof forest $\f$ converges to $0$ as $n\to\infty$. 

\end{corollary}

\proof (i) Let $\kappa$ denote the number of tree-components of $\f$ and 
$f_{1,n}\geq  \cdots \geq f_{\kappa,n}$ the ranked sequence of their sizes. As $U$ and $U'$ are picked uniformly in $[n]$ and independently of $\f$, and $\sum_{i=1}^\kappa f_{i,n}\leq n$, we have
$$\P(\hbox{$U$ and $U'$ are connected in $\f$}) = n^{-2}\E\left(\sum_{i=1}^\kappa f_{i,n}^2\right) \leq \E(n^{-1}f_{1,n})\,.$$
We deduce from Theorem \ref{T2} that $\lim_{n\to\infty}\E(n^{-1}f_{1,n})=1$ in the supercritical case, and therefore
$n^{-1}f_{1,n}$ converges in probability to $1$.

(ii) Clearly there is the lower bound
$$\P(\hbox{$U$ and $U'$ are connected in $\f$}) \geq \varepsilon^2 \P(f_{1,n}\geq \varepsilon n)\,,$$
and we conclude again from Theorem \ref{T2} that the right-hand side converges to $0$ for $\alpha=1/2$. \QED

In the critical case, we point out that despite the absence of giant components, there exists nearly giant tree-components.

\begin{proposition}\label{P1} In the critical case $\alpha=1/2$, for every $\varepsilon >0$, the probability that there exist tree-components of ${\mathcal F}_n$ with size greater than 
$n^{1-\varepsilon}$ converges to $1$ as $n\to\infty$.

\end{proposition}

The rest of this section is devoted first to the proof of Theorem \ref{T2}, and then to that of Proposition \ref{P1}. In this direction, 
 we denote for every $n\in\N$ by $X(n)=X(n,1)$ the number of random cuts which are needed to isolate a typical vertex in a uniform Cayley tree of size $n$,
 using the algorithm of Section 2.1.
Let also $\beta$ stand for a Borel$(1)$ variable wich we assume to be independent of the $X(n)$'s; we shall need an estimate for the Laplace transform of the mixture $X(\beta)$.

\begin{lemma}\label{L4}
 We have
$$\E(1-\exp(-qX(\beta)))\asymp  q \ln(1/q)\,, \qquad q\to 0+,$$
in the sense that the ratio of these two quantities remain bounded away from $0$ and $\infty$ when $q$ tends to $0+$.
\end{lemma}

\proof To ease the notation, $c$ will be used to denote some unimportant numerical constants which are positive and finite, and may be different in different expressions.

First, for every $k\geq 1$, there is the upperbound for the tail of $X(\beta)$
$$\P(X(\beta)>k)\leq \P(\beta>k^2)+ \sum_{n=k+1}^{k^2}\P(\beta=n)\P(X(n)>k)\,.$$
On the one hand, it is well-known from Stirling's formula that  
$$\P(\beta=n)\sim \frac{1}{\sqrt{2 \pi n^3}}\,.$$
 On the other hand, Markov inequality
yields $\P(X(n)>k)\leq k^{-1}\E(X(n))$ and since $\E(X(n))\leq c \sqrt n$ (cf. Meir and Moon \cite{MM}, or Janson \cite{Janson}), we conclude that
$$\P(X(\beta)>k) \leq c k^{-1} \left(1+\sum_{k+1}^{k^2}1/n\right)\leq 2c/k\,.$$
It follows that for $k\geq 2$
$$\E(X(\beta)\wedge k) = \sum_{n=0}^{k-1} \P(X(\beta)>n) \leq c\ln k\,,$$
and therefore for every $0<q\leq 1/2$
$$\E(1-\exp(-qX(\beta)))\leq \E((qX(\beta))\wedge 1) = q\E(X(\beta)\wedge q^{-1})\leq c q \ln (1/q)\,.$$

Second,  there is the lowerbound for the tail of $X(\beta)$
$$\P(X(\beta)>k)\geq \P(\beta>k^2) \, \times\, \inf_{n > k^2}\P(X(n)>k)\,.$$
On the one hand,  $\P(\beta>k^2)\sim \sqrt{2/\pi } k^{-1}$, and on the other hand, since we know from Janson \cite{Janson} and Panholzer \cite{Panholzer} that $X(n)/\sqrt n$ converges weakly to a Rayleigh variable as $n\to\infty$, we also have $\inf_{n > k^2}\P(X(n)>k)\geq c$. Hence 
$$\P(X(\beta)>k)\geq c/k\,,\qquad k\geq 1\,.$$
Now it suffices to write
\begin{eqnarray*}
\E(1-\exp(-qX(\beta))) &=& q \int_0^{\infty} \e^{-qx} \P(X(\beta)>x) \d x\\
&\geq& \e^{-1}q \int_0^{1/q} \P(X(\beta)>x) \d x\,,
\end{eqnarray*}
so the preceding inequality enables us to conclude that
$$\E(1-\exp(-qX(\beta)))\geq c q \ln (1/q)\,,$$
which completes the proof.  \QED

Lemma \ref{L4} yields a rough asymptotic estimate for partial sums of i.i.d. copies of $X(\beta)$ properly conditioned, which will suffice for our purpose.
Specifically, we consider   a sequence $(X_i(n), n\geq 1)_{i\in\N}$ of i.i.d. copies of $(X(n), n\geq 1)$ and a sequence $(\beta_i)_{i\in\N}$ of i.i.d. 
Borel$(1)$ variables which is independent of the preceding. We now state the following technical result.

\begin{corollary}\label{C2} We have for every $0<\varepsilon < 1$ and $a>0$ that:

\noindent {\rm (i)}  if $\alpha >1/2$, then as $n\to\infty$,
 $$
 \max_{ \sqrt{\varepsilon n}  \leq k \leq \sqrt{n/\varepsilon}}
 \P\left( n^{-\alpha}\sum_{i=1}^k X_i(\beta_i) >a \,\Big|\,  \sum_{i=1}^k \beta_i=n  \right)\to 0\,,
 $$

\noindent {\rm (ii)}  if $\alpha =1/2$, then as $n\to\infty$,
 $$
 \max_{ \sqrt{\varepsilon n} \leq k \leq \sqrt{n/\varepsilon}}
 \P\left( n^{-1/2}\sum_{i=1}^k X_i(\beta_i) \leq a \,\Big|\,  \sum_{i=1}^k \beta_i=n  \right)\to 0\,.
 $$

\end{corollary}

\proof  (i) The unconditional version is easy; namely we claim that for $\alpha >1/2$
\begin{equation}\label{E1}
\lim_{n\to\infty} n^{-\alpha}\sum_{i=1}^{\lfloor \sqrt n \rfloor} X_i(\beta_i)=0
 \qquad \hbox{in probability.}
 \end{equation}
Indeed, take  $q>0$ and write
$$
\E\left( \exp\left(-q n^{-\alpha}\sum_{i=1}^{\lfloor \sqrt n \rfloor} X_i(\beta_i)\right)\right)= \left(1-\E\left(1-\exp\left(q n^{-\alpha}X(\beta)\right)\right)\right)^{\lfloor \sqrt n \rfloor}\,.$$
It follows from Lemma \ref{L4}
 that when $n\to \infty$, 
 the quantity above converges to $1$, which establishes \eqref{E1}.
 
 The conditional version now follows from a standard argument based on the local limit theorem. Specifically, write 
 $\tau_k=\beta_1+\ldots + \beta_k$ for the partial sum of i.i.d. Borel variables, so $\tau_k$ has the Borel-Tanner distribution with parameter $k$.
 It is well-known that $k^{-2} \tau_k$ converges in law as $k\to\infty$ to a stable$(1/2)$ random variable. If we denote by $g$ the continuous density of the latter,
 then Gnedenko's  Local  Limit Theorem (see, for instance, Theorem 4.2.1. in \cite{IL}) states that
 \begin{equation}\label{E2}
 \lim_{k\to\infty} \sup_{\ell\geq 1} | k^2\P(\tau_k=\ell)-g(\ell/k^2)| = 0\,.
 \end{equation}
 
 Under the conditional probability given $\tau_k=n$, the $k$-tuple $(X_1(\beta_1), \ldots, X_k(\beta_k))$ is exchangeable, and we have 
 \begin{eqnarray*}
  \P\left( n^{-\alpha}\sum_{i=1}^k X_i(\beta_i) >a \,\Big|\,  \tau_k=n  \right) 
  &\leq & 
  2 \P\left( n^{-\alpha}\sum_{i=1}^{\lceil k/2 \rceil} X_i(\beta_i) >a/2 \,\Big|\,  \tau_k=n  \right)\\
  &\leq & 
  2 \P\left( n^{-\alpha}\sum_{i=1}^{\lceil k/2 \rceil} X_i(\beta_i) >a/2\right)
   \times \sup_{\ell\geq 1}\frac{\P(\tau_{\lfloor k/2\rfloor}=\ell)}{\P(\tau_k=n)}\,.
 \end{eqnarray*}
 On the one hand, we know from \eqref{E1} that the first term in the product above converges to $0$ when $n\to\infty$, uniformly for
 $\sqrt{\varepsilon n}\leq k \leq \sqrt{n/\varepsilon}$. On the other hand, \eqref{E2} yields the bounds
 $$ \max_{ \sqrt{\varepsilon n}  \leq k \leq \sqrt{n/\varepsilon}} \, \sup_{\ell\geq 1} k^2\P(\tau_{\lfloor k/2\rfloor}=\ell) \leq c(\varepsilon)$$
 and (since the stable density $g$ remains bounded away from $0$ on $[\varepsilon, 1/\varepsilon]$)
 $$ \min_{ \sqrt{\varepsilon n}  \leq k \leq \sqrt{n/\varepsilon}}k^2\P(\tau_k=n)\geq c'(\varepsilon)\,,$$
 where $c(\varepsilon)$ and $c'(\varepsilon)$ are two positive constants depending only on $\varepsilon$, which completes the proof.
 
 (ii) The argument is similar for $\alpha=1/2$, except that the unconditional version is now
$$\lim_{n\to\infty} n^{-1/2}\sum_{i=1}^{\lfloor \sqrt n \rfloor} X_i(\beta_i)=\infty\,, \qquad \hbox{in probability.}$$
The conditional version then follows as for (i). \QED

We are now able to proceed to the proof of Theorem \ref{T2}.

{\noindent{\bf Proof of Theorem \ref{T2}:}\hskip10pt}     
Recall the spinal decomposition of the uniform random Cayley tree $\t$  in Section 2.2, and in particular Lemma \ref{L2}.
We consider the dynamics with fires described in the Introduction separately on the spine ${\bf s}$ and each bush
${\bf b}_0, \ldots, {\bf b}_{\lambda_n}$, where $\lambda_n$ is the length of the spine.

Plainly, $U$ and $U'$ are connected in the fireproof forest ${\mathcal F}_n$ if and only if 

$\bullet$ the entire spine ${\bf s}$ is fireproof in the dynamics restricted to ${\bf s}$,

$\bullet$ the root $V_i$ of the bush ${\bf b}_i$ is fireproof in the dynamics restricted to ${\bf b}_i$, for every $i=0, \ldots, \lambda_n$.

We now invoke Lemma \ref{L2} and work conditionally on $\lambda_n=k$ and on the set of vertices ${\mathcal V}_0, \ldots,
{\mathcal V}_k$ of the $k+1$ bushes. Recall that each bush ${\bf b}_i$ rooted at $V_i$ is a uniform rooted Cayley tree on ${\mathcal V}_i$,
independently of the other bushes. In particular, we know from the proof of Theorem \ref{T1} that  
the conditional  probability that $V_i$ is fireproof in the dynamics restricted to
${\bf b}_i$ is
$$\E((1+n^{-\alpha})^{-X(n_i)})\,,$$
where $n_i={\rm Card}({\mathcal V}_i)$ is the size of ${\bf b}_i$ and $X(m)$ stands for the variable that gives the number of random cuts needed to isolate the root
of a uniform rooted Cayley tree of size $m$. Further, the conditional probability that all the $k$ edges of the spine are fireproof
in the dynamics restricted to ${\bf s}$  is 
$(1+n^{-\alpha})^{-k}$. We conclude that the (unconditional) probability that $U$ and $U'$ are connected in ${\mathcal F}_n$ equals
$$\sum_{k=0}^{\infty} \P(\lambda_n=k) (1+n^{-\alpha})^{-k} \E\left(\exp\left(-\ln(1+n^{-\alpha}) \sum_{i=0}^k X_i(\beta_i)\right)\right)\,,$$
where $\beta_i={\rm Card}({\mathcal V}_i)$ is the size of the bush ${\bf b}_i$, and we used the same notation as for Corollary \ref{C2}. 

It is well-known, and also easy to check from Lemma \ref{L2}(i) that $n^{-1/2}\lambda_n$ converges in distribution to a Rayleigh variable as $n\to\infty$.
It now follows from Lemma \ref{L2}(ii) and Corollary \ref{C2} that the quantity above converges as $n\to\infty$ to $1$ if $\alpha >1/2$, and to $0$ if $\alpha = 1/2$. \QED

We now  turn our attention to the proof of Proposition \ref{P1}, which will follow easily from the two following facts.
First, for $1\leq k \leq n$, disconnect a uniform random Cayley tree $\t$ by removing $k-1$ of its edges uniformly at random.
This yields a partition of $[n]$ into $k$ subsets of vertices, say  ${\mathcal V}_1, \ldots, {\mathcal V}_k$, where the labeling is made uniformly at random given the preceding. Denote by $\t^{(i)}$ the restriction of $\t$ to ${\mathcal V}_i$ for $i=1, \ldots, k$.

\begin{lemma}\label{L5}
 \noindent {\rm (i)} Conditionally on the induced partition $[n]= \bigsqcup_{i=1}^k {\mathcal V}_i$,  the $\t^{(i)}$'s
are independent uniform Cayley trees on their respective sets of vertices. 

\noindent {\rm (ii)} The $k$-tuple $(\#{\mathcal V}_1, \ldots,\# {\mathcal V}_k)$ of the sizes of the induced partition is distributed as $k$ independent Borel$(1)$ variables conditioned to have sum $n$.
\end{lemma}

\proof For $k=2$, the first claim reduces to Lemma \ref{L2}(iii), and the general case follows by induction. The second assertion is known from Pavlov \cite{Pavlov} 
and Pitman \cite{Pi}; see also Corollary 5.8 in \cite{RFC}. \QED

We next recall informally a well-known property of increasing random walks  with a step distribution which has a regularly varying tail of exponent $>-1$. 
Conditioning such a random walk to be abnormally large after $k$ steps where $k$ is large, essentially amounts to conditioning the random walk
on having one of its steps abnormally large. Here is a consequence tailored for our need.

\begin{lemma}\label{L6} Let $\beta_1, \ldots$ be a sequence of i.i.d. Borel$(1)$ variables and set $\tau_k=\beta_1+\cdots + \beta_k$.
Fix $\varepsilon >0$; let  $k(n)=\lfloor n^{(1-\varepsilon)/2}\rfloor$ and denote by $\beta^*_2$ the second largest value amongst $\beta_1, \ldots, \beta_{k(n)}$. Then
$$\lim_{n\to\infty}\P(\beta^*_2\in [n^{1-2\varepsilon}, n^{1-\varepsilon/2}] \mid \tau_{k(n)}=n)=1\,.$$
\end{lemma}

\proof We first recall that $\tau_{k(n)}$ has the Borel-Tanner distribution with parameter $k(n)$, so that an application of Stirling formula
readily yields
\begin{eqnarray*}
\P(\tau_{k(n)}=n)=\frac{k(n)}{(n-k(n))!} \e^{-n} n^{n-k(n)-1} &\sim & (2\pi n^{2+\varepsilon })^{-1/2}\\
&\sim& k(n) \P(\beta=n)\,.
\end{eqnarray*}

Next, note that
$$\ln \P(\max_{i=1, \ldots , k(n)-1}\beta_i \leq n^{1-2\varepsilon}) \sim (k(n)-1) \ln\left(1-(2\pi n^{1-2\varepsilon})^{-1/2}\right)
\sim -n^{\varepsilon/2}\,,$$
which entails
$$\P(\beta^*_2\leq n^{1-2\varepsilon})\leq k(n)  \P(\max_{i=1, \ldots , k(n)-1}\beta_i \leq n^{1-2\varepsilon}) =o(\P(\tau_{k(n)}=n))\,.$$
As a consequence 
$$\lim_{n\to\infty} \P(\beta^*_2 \leq n^{1-2\varepsilon}\mid \tau_{k(n)}=n)=0\,.$$

Then we write
\begin{eqnarray*}
& & \P(\max_{i=1, \ldots , k(n)-1}\beta_i \leq n^{1-\varepsilon/2}, \tau_{k(n)-1}\leq k^2(n)\ln n , \tau_{k(n)}=n) \\
&=&\sum_{j=k(n)-1}^{k^2(n)\ln n}  \P(\max_{i=1, \ldots , k(n)-1}\beta_i \leq n^{1-\varepsilon/2}, \tau_{k(n)-1}=j)
\P(\beta=n-j)\\
&\sim& \frac{1}{\sqrt{2\pi n^3}} \P(\max_{i=1, \ldots , k(n)-1}\beta_i \leq n^{1-\varepsilon/2}, \tau_{k(n)-1}\leq k^2(n)\ln n)\,,
\end{eqnarray*}
where at the second line we used that 
$$\P(\beta=n-j)\sim (2\pi n^3)^{-1/2}\hbox{  uniformly for }k(n)-1\leq j \leq k^2(n)\ln n\,.$$
It is immediately checked that 
$$ \P(\max_{i=1, \ldots , k(n)-1}\beta_i \leq n^{1-\varepsilon/2})\sim 1\,,$$
and, because $k(n)^{-2}\tau_{k(n)-1}$ converges in distribution to a stable$(1/2)$ variable, we also have
$\P(\tau_{k(n)-1}\leq k^2(n)\ln n)\sim 1$. Putting the pieces together, we get
$$\P(\max_{i=1, \ldots , k(n)-1}\beta_i \leq n^{1-\varepsilon/2}, \tau_{k(n)-1}\leq k^2(n)\ln n , \tau_{k(n)}=n)\sim \frac{1}{\sqrt{2\pi n^3}}\,.$$

We conclude from an argument of symmetry that
$$
\P(\beta^*_2\leq n^{1-\varepsilon/2}, \tau_{k(n)}=n)
\geq  k(n) \P(\max_{i=1, \ldots , k(n)-1}\beta_i \leq n^{1-\varepsilon/2}, \tau_{k(n)-1}\leq k^2(n)\ln n , \tau_{k(n)}=n)$$
and then
$$\liminf_{n\to\infty} 
\P(\beta^*_2\leq n^{1-\varepsilon/2} \mid \tau_{k(n)}=n) 
\geq \liminf_{n\to \infty} \frac{k(n)}{\sqrt{2\pi n^3} \P(\tau_{k(n)}=n)} \sim 1\,;
$$
which completes the proof. \QED

We are now all the ingredients to establish Proposition \ref{P1}.

{\noindent{\bf Proof of Proposition \ref{P1}:}\hskip10pt}    We shall establish the claim with $n^{1-3\varepsilon}$ replacing $n^{1-\varepsilon}$ in the statement,
which makes no difference as $\varepsilon>0$ can be chosen arbitrarily small.

Recall that $\alpha=1/2$ and take the same notation as in Lemma \ref{L6}.
We start by considering the dynamics with fires described in the Introduction until the instant when the $k(n)$-th inflammable edge of $\t$ is replaced by a fireproof one.
Because $k(n)\ll \sqrt n$, the probability that a fire has occurred before that instant is small when $n$ is large, and we implicitly 
work conditionally on the complementary event in the sequel. 

Let us now remove these $k(n)$ first fireproof edges, and focus on ${\mathcal V^*_2}$, the second largest set of vertices in the induced partition.
According to Lemma \ref{L5}(i), conditionally on $\#{\mathcal V^*_2}=m$, the restriction of $\t$ to ${\mathcal V^*_2}$ is a uniform Cayley tree on a set of $m$ vertices. 

Because the edge connecting ${\mathcal V^*_2}$ to $\t \backslash {\mathcal V^*_2}$ is fireproof, no fire started outside ${\mathcal V^*_2}$
can burn edges between vertices of ${\mathcal V^*_2}$. 
We also make the following elementary observation of monotonicity. One can couple two fire-dynamics on the same tree with different 
firing rates such that every edge which is burnt for the dynamics with the smallest firing rate is also burnt for the dynamics with the highest firing rate.
As a consequence, a tree-component in the fireproof forest for the dynamics with the highest firing rate is always contained into a tree-component
of the fireproof forest for the dynamics with the smallest firing rate.

Thanks to Lemmas \ref{L5}(ii) and \ref{L6}, we may focus on the case when the size of ${\mathcal V^*_2}$ belongs to  $[n^{1-2\varepsilon},  n^{1-\varepsilon/2}]$. Now conditionally on $\#{\mathcal V^*_2}=m$, we  run the fire dynamics restricted to ${\mathcal V^*_2}$
with firing rate $m^{-1/(2-\varepsilon)}\geq n^{-1/2}$. Since $1/(2-\varepsilon)>1/2$, we deduce from Corollary \ref{C1}(i) and the monotonicity property observed above
that the conditional probability
given $\#{\mathcal V^*_2}$
that the fireproof forest ${\mathcal F}_n$ restricted to  ${\mathcal V^*_2}$ contains a tree-component of size at least $\frac{1}{2} \#{\mathcal V^*_2}$, 
converges in probability to $1$. A fortiori,  the conditional probability
given $\#{\mathcal V^*_2}$
that the fireproof forest restricted to  ${\mathcal V^*_2}$ contains a tree-component of size at least $n^{1-3\varepsilon}$, 
converges in probability to $1$, which completes the proof of our statement. \QED
\end{section}

\vskip 5mm

\noindent{\bf Acknowledgments}.  
This work has been  realized during a stay at the
 Forschungsinstitut f\"ur Mathematik, ETH Z\"urich, 
whose support is thankfully acknowledged.
It is also supported by ANR-08-BLAN-0220-01.

  \end{document}